\documentclass[12pt]{amsart}
\usepackage{subfigure}
\usepackage{graphicx}
\usepackage{rotating}
\usepackage{morefloats}
\usepackage[utf8]{inputenc}
\usepackage[T1]{fontenc}
\usepackage{amssymb}
\usepackage{amsmath}    
\def\R{I\kern -0,37 em R}
\def\P{I\kern -0,37 em P}
\def\Z{I\kern -0,37 em Z}

\begin{document}

\title[Pfaffian Systems]{Non-Integrable Pfaffian Systems}
\author{A. Kumpera}
\address[Antonio Kumpera]{Campinas State University, Campinas, SP, Brazil}
\email{antoniokumpera@hotmail.com}

\date{July, 2016}


\maketitle

\begin{abstract}
We discuss a recurrent geometrical method, due to Élie Cartan and von Weber (\cite{Cartan1901},\cite{Weber1898}) enabling us to determine, step by step, the maximal  integral manifolds of a not necessarily integrable nor regular  Pfaffian system. The dimensions of such integral manifolds can, of course, vary from point to point but more so can vary at a given point it depending upon the choice of their recurrent buildup. When the system is regular and integrable then, of course, we obtain the maximal integral leaves of the integral foliation. Attention is also given to those integrable systems that can be integrated by quadratures which was, in the 19th century, the dream of many. However, our main interest resides in enhancing the \textit{Jordan-Hölder} integration procedure so as to construct the local maximal integral manifolds that find many applications.
\end{abstract}

\section{Introduction}
The initial purpose, of the present discussion is to search for point-wise criteria that will detect the maximal dimensions of the integral manifolds passing through given points, the Pfaffian systems being generic \textit{i.e.}, not necessarily integrable (involutive). Subsequently, we make use of the \textit{Jordan-Hölder} integration procedure (\cite{Kumpera2016}) so as to construct or at least provide a convenient description for the local maximal integral manifolds. This procedure can already be found, though veiled in behind the curtains, in most of Cartan's later work. We also hope that, in a forthcoming paper, attention be given to global integral manifolds. Curiously enough, globalization makes appeal to topology and seems far more delicate and involved than geometry. Given a Pfaffian system $\mathcal{P}$, we start building, at a given though arbitrary point $x_0$, \textit{integral} linear contact elements thus obtaining an ascending chain of sub-spaces that forcibly ends up with an element that cannot be further enlarged, the chain starting with the null space followed by a \textit{linear} ($1-$dimensional) integral element and so forth. We finally observe that a substantial part of our discussion and mainly that concerning the integral elements can also be found, partially, in \cite{Cartan1900}. 

\section{Integral contact elements}
Let \textit{M} be a finite dimensional differentiable manifold and $\mathcal{P}$ a Pfaffian system of constant rank \textit{r} defined on \textit{M}, assumed to be of dimension \textit{n}. For convenience, we also assume that all the data are of class $\mathcal{C}^{\infty}$ though, of course, much weaker differentiability requirements could be assumed. We further denote by $\{\omega^1,\omega^2,~\cdots~,\omega^r\}$ a free system of local generators, defined in a neighborhood of a generic point $x_0$, and by $\Sigma$ the annihilator of $\mathcal{P}$ in $TM,$ thus obtaining the field of $(n-r)-$dimensional first order contact elements annihilated by $\mathcal{P}.$

\vspace{5 mm}
\noindent
Let us take an arbitrary point $x_0\in M$ and try to determine, \textit{de proche en proche} as would say Cartan, an ascending chain of integral contact elements at the point $x_0.$ We firstly take a one-dimensional subspace $l\subset \Sigma_{x_0},$ a line in the terminology of Cartan, and examine when such a linear subspace can be included in a higher dimensional \textit{integral} contact element. The sole condition $l\subset\Sigma_{x_0}$ is not sufficient and the entire subsequent discussion reposes on the following remark:

\vspace{5 mm}
\noindent
\textit{Given an integral sub-manifold $\mathcal{I}$ of the Pfaffian system $\mathcal{P}$ and an arbitrary local section $\omega$ of this system, the restrictions $\iota^*\omega$ and $\iota^*d\omega$ vanish identically on $\mathcal{I},$ where $\iota:\mathcal{I}\longrightarrow M$ is the inclusion.}

\vspace{5 mm}
\noindent
It suffices, of course, to check the above condition for the generators $\omega^i$ and, further, this condition transcribes by the two requirements $i(v)\omega^i=0,~v\in l,~v\neq 0,$ and $i(v)d\omega^i\equiv 0\hspace{3 mm}mod~\{\omega^i\},$ the latter condition also rewriting by $(d\omega^i)_{x_0}\equiv 0\hspace{3 mm}mod~\{\omega^i\}.$

\vspace{5 mm}
\noindent
We now fix a given $1-dimensional$ linear integral contact element $E_1$ issued at the point $x_0$ (\textit{i.e.}, contained in $\Sigma_{x_0}$) and look for the conditions under which a plane contact element $E_2$ ($dim=2$) containing the given linear element $E_1$ is also integral. In order to find such an element, it suffices to look for a vector $w\in T_{x_0}$ linearly  independent from \textit{v} and verifying the two conditions $i(w)\omega^i=0,~ i(w)i(v)d\omega^i=d\omega(v,w)=0$. We shall say that \textit{w} is \textit{in involution} with \textit{v} (in Cartan's terminology, the line $\lambda$ generated by \textit{w} is in involution with the line \textit{l}). Linearity then implies that all the vectors \textit{w} that are in involution with \textit{v} (adding, of course, the null vector) form a linear subspace $\tilde{E}_1$ containing all the $2-$dimensional integral subspaces considered above. 

\vspace{5 mm}
\noindent
Let us next fix a $2-dimensional$ integral contact element $E_2$ (a linear $2-$dimensional subspace of $\tilde{E}_1$) and consider the vectors
\textit{z} that are in involution with all the vectors belonging to $E_2.$ It suffices, of course, to consider those vectors that are in involution, simultaneously, with those belonging to any given basis of $E_2$ and we shall denote by $\tilde{E}_2$ the linear subspace consisting of all these vectors that obviously is contained in $\tilde{E}_1.$

\vspace{5 mm}
\noindent
Fixing a $3-dimensional$ contact element $E_3$ that is contained in  $\tilde{E}_2$ and contains $E_2,$ we can pursue this nice \textit{Origami} game and construct an ascending chain $E_0\subset E_1\subset E_2\subset E_3~ \cdots~,$ each element satisfying the involutiveness condition with respect to the preceding one. We include the null space $E_0$ just to please Cartan where it figures, in his writings, as the one point set $E_0=\{x_0\}.$ Dimensions being finite, there will be an integer $\rho_k$ for which the above process breaks down namely, there will not exist any other non-null vector, in involution with $E_{\rho_k},$ other than those already obtained previously. 

\newtheorem{uvw}[DefinitionCounter]{Definition}
\begin{uvw}
 According to Cartan, we shall say that the integer $s_1=n-r-\rho_k$ is the \textbf{character} (caractère) of $\mathcal{P}$ at the point $x_0$ relative to the above ascending chain.
\end{uvw}

\noindent
Described more geometrically, the character is simply the difference between the values of $dim~\Sigma_{x_0}$ and $dim~E_{\rho_k}$. We next describe Cartan's most ingenious method for determining these maximal integral elements. Let us take a tangent contact element \textit{E} contained in some $\Sigma_{x_0}$ as well as a second sub-space $E'\subset E.$ The sub-space $E'$ is said to be \textit{characteristic} with respect to $E$ whenever any element of $E'$ is in involution with every element of $E.$ The sum of all the characteristic sub-spaces is then a maximal characteristic sub-space called the \textit{characteristic element} of $E$ with respect to $\mathcal{P}.$ Though $\mathcal{P}$ is assumed to have null characteristics, these characteristic elements need not be null. Next, given any two sub-spaces \textit{E} and \textit{F} contained in $\Sigma_{x_0},$ we shall say that they are \textbf{conjugate} when their characteristic elements are in involution, this meaning that any vector belonging to one of them is in involution with all the vectors belonging to the other. In order to obtain an integral element of dimension $\mu\leq\rho_k,$ it suffices to take $\mu$ independent vectors $\{v_1,~\cdots~,v_{\mu}\}$ for which the corresponding sub-spaces $\tilde{E}_{2,i},~1\leq i\leq\mu,$ are pairwise conjugate and then simply take their intersection.

\vspace{5 mm}
\noindent
We next define $s_j=dim~\tilde{E}_j,$ for  $1\leq j\leq\rho_k$ and at the given point $x_0$. Obviously, $s_1-s_2\geq s_2-s_3\geq~\cdots~\geq s_{\rho_k-1}-s_{\rho_k}$ since, at each step, we have a lesser choice of vectors to look for and, furthermore, $\rho_k=n-r$ if and only if $\mathcal{P}$ is integrable in a neighborhood of $x_0.$ We should also observe that this point-wise technique is most appropriate for the study of singular Pfaffian systems in view of obtaining the maximal integral sub-manifolds whose dimensions can vary from point to point. The integers $s_j$ are called the successive \textit{enlarged characters} of $\mathcal{P}$ at the point $x_0.$ We finally call the attention of the reader to the way Cartan embraced his thoughts. Invariably and persistently, he would always pursue the \textit{co-variant} trail and, in his \textit{mémoire} \cite{Cartan1901} as well as previously in \cite{Cartan1900}, this shows up as soon as he describes the above-mentioned contact elements. In fact, these are depicted with the help of convenient linear combinations of the differentials $dx^i$ relative to a local coordinate system $(x^1,~\cdots~,x^n)$ defined in a neighborhood of the point $x_0.$ In other terms, Cartan employs linear bases for the sub-spaces, in $T^*M,$ that annihilate the \textit{contra-variant} subspaces $E_j$ and $\tilde{E}_j$ considered previously ([\textit{nous}] \textit{regardons $dx_1,dx_2,~\cdots~,dx_r$ comme les paramètres directeurs d'une droite issue de ce point}, pg.4, \textit{l.}8). In so doing, the calculations become much simpler and completely straightforward though, unfortunately, our \textit{contra-variant} intuition of the geometric world around us suffers with it.

\section{Integral sub-manifolds}
We first observe that any Pfaffian system whose rank does not exceed $n-1$ always admits integral $1-$dimensional sub-manifolds, at a given point $x_0,$ since it suffices to integrate any local vector field $\xi$ that is annihilated by the system and that satisfies $\xi_{x_0}\neq 0.$ The above assertion can fail for singular systems for, in this case, it might happen that $\Sigma_{x_0}=0$ at a given point $x_0.$ Hence, there is in principle no integrability condition pending for the dimension 1. Secondly, we observe that the dimensions of all the previously considered contact elements are \textit{lower semi-continuous} (the dimensions tend to increase) and, further, that the \textit{pertinence relation} $\in$ is a \textit{closed} relation. Consequently, these properties enable us to choose locally, in a neighborhood of $x_0,$ a finite set of linearly independent vector fields $\xi_{\mu},~1\leq\mu\leq\rho_k,$ such that

\vspace{5 mm}
\textbf{1.} The vectors $\{\xi_{\mu}(x_0)\}$ generate $E_{\rho_k}$ and

\vspace{3 mm}
\textbf{2.} The vector fields $\xi_{\mu}$ generate, in that neighborhood, an integrable distribution (field of contact elements).

\vspace{5 mm}
\noindent
We can then integrate this local distribution and obtain an integral sub-manifold $\mathcal{I}$ of the initially given Pfaffian system $\mathcal{P}$, this integral manifold being of maximal dimension relative to the point $x_0$ and to the ascending chain. We should however observe that it is not unique relative to the point $x_0$ unless, of course, its dimension be equal to the co-rank of $\mathcal{P}.$ When this is not the case, there are many possible choices for the chains. Since the element $E_{\rho_k}$ is maximal, relative to the sub-chain terminating by $E_{\rho_k-1},$ this element is unique, always relative to the sub-chain, though the above indicated integral sub-manifold needs not be locally unique, relative to the complete chain. Furthermore, choosing the local vector fields $\xi_{\mu}$ in such a way that the vectors $\{\xi_{\mu}(x_0),~\mu\leq j\}$ generate $E_j,$ we can obtain, by successive integrations, an ascending chain of integral sub-manifolds respectively tangent to the elements of the linear chain, each one contained in the next manifold. We finally observe that the singularities as well as the singular integral sub-manifolds can be located inasmuch as evaluated by computing the integer $\rho_k$ at the various points of the manifold \textit{M} and corresponding to the various ascending chains of integral elements.

\section{Characteristics and derived systems}
Given a Pfaffian system $\mathcal{P},$ we recall that the \textit{characteristic system} associated to $\mathcal{P}$ is the Pfaffian system generated by all the differential 1-forms $\omega$ and $i(\xi)d\omega$ where $\xi$ is an arbitrary local vector field annihilated by $\mathcal{P}$ and $\omega$ an arbitrary local section of the given system. The contra-variant counterpart namely, the annihilator of the characteristic system is defined as follows: At each point $x\in M,$ we consider the sub-space, of the tangent space, whose elements are the vectors \textit{v} such that $i(v)\omega=0$ and $i(w)i(v)d\omega=0$ for any vector \textit{w} annihilated by $\mathcal{P}$ and any local $1-$form $\omega$ belonging to the system. In other terms, the above \textit{characteristic} sub-space is defined as the set of all the vectors \textit{v} annihilated by $\mathcal{P}$ that further verify $d\omega(v,w)=0$ for any vector \textit{w} also annihilated by $\mathcal{P}$ and any local section $\omega$ of $\mathcal{P}.$ We denote by $\mathcal{CH}$ this characteristic system and observe that its purpose resides in (loosely speaking) finding the minimum number of independent function by means of which the initially given system can be re-written or, equivalently, in finding the local fibration with the largest dimensional fibres with respect to which the given system factors to the quotient \textit{i.e.,} to the base space of this fibration onto an equivalent \textit{quotient} system (\textit{cf.} \cite{Kumpera1982}). Inasmuch, the \textit{derived system} of $\mathcal{P},$ denoted by $\mathcal{P}_1,$ is the Pfaffian system generated by the forms $\omega,$ local sections of $\mathcal{P},$ for which $d\omega\equiv 0\hspace{3 mm}mod~\mathcal{P}$ \textit{i.e.,} generated by those local forms that satisfy the integrability condition with respect to $\mathcal{P}.$ The derived system is, of course, a Pfaffian sub-system of the initially given system. Iterating this procedure, we obtain a descending chain of Pfaffian systems terminating by an integrable system (its derived system coincides with the system) that, eventually, reduces to the null system. The contra-variant description of the derived system is left to the reader though, in the author's opinion, it is rather clumsy to deal with (\textit{cf.} \cite{Kumpera2014}, \cite{Monty2001} and \cite{Mormul2005}).

\vspace{5 mm}
\noindent
We can now apply the above inquiry, concerning integrable elements, to both characteristic and derived systems though, presently, we shall only highlight the most relevant facts and, among them, those facts that were placed in evidence by Cartan (pg.9, n.6 and the following pages). We first observe, by simply inspecting the contra-variant definition, that the characteristic system is equal, at each point, to the union of those  $1-$integral elements (or tangent vectors) of  $\mathcal{P}$ that are in involution with all the other $1-$integral elements of $\mathcal{P}.$ The characteristic system is always integrable (involutive) independently of the nature of $\mathcal{P}$ and, if moreover it is regular, we can factor locally its integral foliation. Furthermore, we can also factor all the differential forms of $\mathcal{P},$ \textit{modulo} the leaves, and thereafter the initially given system factors as well to an equivalent quotient system  $\tilde{\mathcal{P}}.$ The integral manifolds of the initial systems are the pullbacks (inverse images) of the integral manifolds in the quotient. Let us also observe that the quotient system has null characteristics.

\vspace{5 mm}
\noindent
We next give a glance at the derived system. From a \textit{contra-variant} point of view, the contact element determined by $\mathcal{P}_1$ at the point $x_0$ is the linear subspace generated by all the vectors annihilated by $\mathcal{P}$ together with all the brackets $[\xi,\eta]_{x_0}$ where $\xi$ and $\eta$ are two local vector fields annihilated by $\mathcal{P}$. In the \textit{contra-variant} mode, the chain of successive derived spaces or systems is increasing whereas, \textit{co-variantly}, it decreases.

\vspace{5 mm}
\noindent
Obviously,
\begin{equation*}
\mathcal{P}_1\subset\mathcal{P}\subset\mathcal{CH}
\end{equation*}
hence, in what concerns the corresponding pseudo-groups of local automorphisms and taking into account the functoriality of the constructions, we obtain the following two inclusions: 
\begin{equation}
\Gamma_{\mathcal{P}}\subset\Gamma_{\mathcal{CH}}\hspace{10 mm}and\hspace{10 mm}\Gamma_{\mathcal{P}}\subset\Gamma_{\mathcal{P}_1}~.
\end{equation}
We can then proceed to extend the above second inclusion so as to obtain a composition series and, eventually, a \textit{Jordan-Hölder resolution}. Making use of the results in \cite{Kumpera2016}, we might then simplify the integration process of $\mathcal{P}.$

\section{The pseudo-groups of local automorphisms}
The hypotheses and the notations being as above, we shall now have a closer look into the integration process of the Pfaffian system $\mathcal{P}.$ According to Lie and Cartan, we first have to choose convenient intermediate pseudo-groups, in the present case these will be associated to structures, in such a way that the second inclusion in (1) extends hopefully to a Jordan-Hölder resolution. As a first step, we examine the local automorphisms of the successive derived Pfaffian systems and compare these pseudo-groups with each other as well as with the automorphisms of the initial system $\mathcal{P}.$

\vspace{5 mm}
\noindent
Let $\mathcal{P}$ be a Pfaffian system defined on the $n-$dimensional manifold \textit{M} and let us assume, at least for the time being, that it is regular \textit{i.e.,} that $rank~\mathcal{P}_x=rank~\mathcal{P}_y$ for any two arbitrary points $x,y\in M$ where, in the \textit{co-variant} jargon, the \textit{rank} is equal to the \textit{dimension}. We then consider a local diffeomorphism (referred to as a local transformation) $\varphi:\mathcal{U}\longrightarrow\mathcal{V},$ where $\mathcal{U}$ and $\mathcal{V}$ are open subsets in \textit{M} and the transformation $\varphi$ is assumed to be differentiable. To say that $\varphi$ preserves or leaves invariant the Pfaffian system $\mathcal{P}$ means, of course, that $\varphi^*(\mathcal{P}_y)=\mathcal{P}_x$ as soon as $\varphi(x)=y.$ We shall also refer to $\varphi$ as being an \textit{automorphism} of $\mathcal{P}.$ Observing that the linear automorphism
\begin{equation*}
\varphi^*_x:\mathcal{P}_y\longrightarrow\mathcal{P}_x
\end{equation*}
is entirely determined by its $1-$jets $j_1\varphi(x)$ (in fact, it is essentially the same thing), we infer that any local transformation $\varphi,$ for which $j_1\varphi(x)$ transposes to an automorphism pertaining to $L(\mathcal{P}_y,\mathcal{P}_x),$ belongs to $\Gamma_{\mathcal{P}}$ and,  consequently, this pseudo-group is a Lie pseudo-group of order 1.\footnote{The groupoid of first order jets is a differentiable manifold, the source and target maps are submersions (more precisely, surmersions) and the composition operation as well as the passage to inverses are both differentiable mappings the latter being a diffeomorphism.}

\vspace{5 mm}
\noindent
Our next task is to show that the Lie groupoid composed by all the $1-$jets of the elements belonging to $\Gamma_{\mathcal{P}}$ is an involutive first order partial differential equation in the sense of Élie Cartan (\cite{Kumpera1972}). Let us denote by $\textbf{G}^1_{\mathcal{P}}$ this first order groupoid associated to $\Gamma_{\mathcal{P}}$ and let us recall that the \textit{standard prolongation} of the differential equation $\textbf{G}^1_{\mathcal{P}}$ is the set of all second order jets $j_2\psi(x),$ of local maps $\psi$ of \textit{M}, such that the \textit{Ehresmann flow}
\begin{equation*}
y\longmapsto j_1\psi(y)
\end{equation*}
is tangent, at first order and at the point \textit{x}, to the equation $\textbf{G}^1_{\mathcal{P}}.$ This means, in denoting by $\tilde{\psi}$ the above flow, that the image sub-manifold of $\tilde{\psi}$ is tangent to $\textbf{G}^1_{\mathcal{P}}$ at the point $\tilde{\psi}(x)$ ($\tilde{\psi}_*(T_xM)\subset T_{\tilde{\psi}(x)}\textbf{G}^1_{\mathcal{P}}$). We can define inasmuch the standard prolongations of the higher order jet groupoids associated to the pseudo-group $\Gamma_{\mathcal{P}}$ and a straightforward argument shows the following lemma.

\vspace{2 mm}
\newtheorem{jhi}[LemmaCounter]{Lemma}
\begin{jhi}
All the $k-$th order groupoids, $k>1,$ associated to the pseudo-group $\Gamma_{\mathcal{P}}$ are the successive standard prolongations of the first order groupoid $\textbf{G}^1_{\mathcal{P}}$
\end{jhi}

\vspace{5 mm}
\noindent
More important, we want to show that the first order groupoid as well as all the higher order groupoids are \textit{involutive equations} in the sense of Cartan. The above Lemma being a first step, we now have to show that all the linear \textit{symbols} of these equations are involutive \textit{i.e.,} that they are $2-acyclic.$ A rather long juggling with the techniques found in \cite{Kumpera1972}  including the application, to the previous equations, of the \textit{Spencer linear complex} will eventually prove the following result.\footnote{We do not provide the details of the proof since it is rather long and the result, in itself, is just auxiliary and not of first order importance.}

\vspace{2 mm}
\newtheorem{klm}[LemmaCounter]{Lemma}
\begin{klm}
All the groupoids associated to the Lie pseudo-group  $\Gamma_{\mathcal{P}}$ are involutive.
\end{klm}

\vspace{5 mm}
\noindent
Let us next consider the derived sequence associated to $\mathcal{P}.$ It is a descending chain
\begin{equation*}
\mathcal{P}\supset\mathcal{P}_1\supset~\cdots~\supset\mathcal{P}_{\mu}=\mathcal{P}_{\mu+1}
\end{equation*}
the last term being integrable or eventually null and, by the naturality of the constructions, we conclude that
\begin{equation}
\Gamma(\mathcal{P}_{\mu})\supset~\cdots~\supset\Gamma(\mathcal{P}_1)\supset\Gamma(\mathcal{P})
\end{equation}
and consequently that
\begin{equation}
\textbf{G}^k_{\mathcal{P}_{\mu}}\supset~\cdots~\supset\textbf{G}^k_{\mathcal{P}_1}\supset\textbf{G}^k_{\mathcal{P}}
\end{equation}
for all integers $k\geq 1.$ In general, there is no reason for the previous sequence to be a composition series (each term being normal in the preceding one) though we shall exhibit a rather outstanding situation where not only this is true but, furthermore, the series becomes \textit{Jordan-Hölder}. 

\vspace{5 mm}
\noindent
Let us now assume that the Pfaffian system $\mathcal{P}$ is a \textit{flag system} \textit{i.e.,} $rank~\mathcal{P}_i=rank~\mathcal{P}_{i+1}+1$ for all \textit{i}. Then, with the notations of the reference \cite{Kumpera2014}, the finite and infinitesimal automorphisms of any of the factored derived systems  $\overline{\mathcal{P}}_{\mu}$ (\textit{modulo} their characteristic variables) is canonically equivalent (isomorphic), via a natural \textit{merihedric} prolongation algorithm, to those of $\mathcal{P}$ and consequently all the corresponding higher order groupoids are canonically isomorphic. However, when these derived systems are considered on their \textit{ambient} space \textit{M}, new automorphisms do appear namely, those of the associated characteristic systems that are also automorphisms of the derived systems in evidence. However, a characteristic system being always integrable, its automorphisms are, locally, those that leave invariant (\textit{i.e.,} permute) the fibres of the fibration obtained by integrating the characteristics and, still more locally, those that leave invariant (permute) the parallel spaces to one of the components in a product. In terms of local coordinates, we consider all the local transformations in an $n-$space, a \textit{simple} pseudo-group according to Cartan, and just retain those transformations that maintain a certain number of coordinates depending only upon these restricted  coordinates, the resulting pseudo-group being \textit{simple} as well. The quotients of these pseudo-groups and, inasmuch, the quotients of the groupoids are therefore simple and the sequence
\begin{equation}
\textbf{G}^k_{\mathcal{P}_{\mu}}\supset~\cdots~\supset\textbf{G}^k_{\mathcal{P}_1}\supset\textbf{G}^k_{\mathcal{P}}
\end{equation}
is Jordan-Hölder. The same argumentation holds on the infinitesimal level and the integration (or not) of the system $\mathcal{P}$ will entirely depend upon the nature of the quotient groupoids and the algebroids.

\vspace{5 mm}
\noindent
Let $E_{\rho_k}$ be a maximal integral contact element of $\mathcal{P}$ at the point $x_0$ and $\mathcal{I}_{\rho_k}$ a local integral sub-manifold of maximal dimension, tangent to $E_{\rho_k}.$ Then, of course, $\mathcal{I}_{\rho_k}$ is also an integral sub-manifold of the derived system, though not forcibly maximal, and will be contained in maximal integral manifolds (not all) of this derived system.

\section{Pfaffian systems whose characters are equal to one}
We now assume that $\rho_k=n-r-1,$ where we recall that $n=dim~M$ and $r=rank~\mathcal{P}.$ In other terms, the dimension of $E_{\rho_k}$ and consequently also that of $\mathcal{I}_{\rho_k}$ are one unit less than the dimension of $\Sigma_{x_0},$ the contact element at the point $x_0$ annihilated by $\mathcal{P}_{x_0}.$ Let us show that $rank~\mathcal{P}_1=r-1$ and provide a criterion for the integrability of this derived system. Cartan and, earlier, von Weber already gave such a criterion that Cartan claimed to be \textit{un fait très remarquable} (pg.25, \textit{l.}-4, Théorème).  

\vspace{5 mm}
\noindent
Firstly, since $\mathcal{P}$ is not integrable (involutive) on account of the maximality of $E_{\rho_k},$ it follows that $rank~\mathcal{P}_1\leq r-1.$ On the other hand, since the dimension of $\mathcal{I}_{\rho_k}$ is equal to $n-r-1$ and all the forms $\omega^i$ vanish on this integral manifold, we can assume relabeling if necessary by $\{\varphi^i\}$ the local basis of $\mathcal{P},$ that $\mathcal{P}_1$ is generated by $\{\varphi^2,~\cdots~,\varphi^r\}.$ In fact, choosing a local section $\varphi_1$ of $\mathcal{P}$ verifying $d\varphi_1|\Sigma_{x_0}\neq 0,$ we can argue together with Cartan (pg.20,21) without performing the calculations.

\vspace{5 mm}
\noindent
Secondly and returning to the consideration of the \textit{integral} contact elements, we consider at the point $x_0$ a maximal integral contact element $E_{\rho_k}.$ As evidenced previously, this element is equal to the tangent space, at $x_0,$ of an integral sub-manifold \textit{W} of the system $\mathcal{P}.$ Locally, this sub-manifold is integral with respect to the forms $\omega^i,$ of a local basis of $\mathcal{P},$ and some additional independent forms $\varphi^s$ that restrict to a single form $\varphi$ in the case when the \textit{character} of the system is equal to one. But then, for the Pfaffian system $\mathcal{Q}$ generated by these forms, we can write
\begin{equation*}
d\omega^i\equiv c^i\mu^i\wedge\varphi\hspace{5 mm}mod~\mathcal{P},
\end{equation*}
with certain functions $c^i$. Consequently, the forms belonging to the sub-module of $\Gamma(\mathcal{P})$ defined by the equation $~\sum~c^i\omega^i=0~$ are precisely those that belong to $\mathcal{P}_1.$ This system is under certain conditions integrable namely, when its \textit{gender} is larger than 1. This invariant, already introduced by von Weber, has to do with the "number" of differential $1-$forms independent from those forming $\mathcal{P}$ and that enter in the expressions of $d\omega\hspace{2 mm}mod~\mathcal{P},~\omega\in\Gamma{\mathcal{P}}$ (local sections). We shall consider this invariant only further but define it right below. Taking a maximal integral manifold of $\mathcal{P}_1,$ when this system is integrable, and containing $x_0,$ we can restrict all the data to this manifold, the investigation of the properties of $\mathcal{P}$ becoming considerably facilitated. We shall also examine further those Pfaffian systems with characters equal to two or more but presently we want to exhibit Cartan's \textit{co-variant} version of the former discussion.

\vspace{2 mm}
\newtheorem{abc}[DefinitionCounter]{Definition}
\begin{abc}
The \textbf{gender} (genre) of a local exterior differential form $\Omega,$ with respect to a Pfaffian system $\mathcal{P},$ is the smallest integer \textbf{h} such that $\Omega^{h+1}\equiv 0\hspace{2 mm}mod~\mathcal{P}$ (the exponent refers to wedge products). The \textbf{gender} of a family of local exterior differential forms is the maximum value of the gender of its elements. 
\end{abc}

\vspace{2 mm}
\noindent
When $\mathcal{P}$ is locally generated by $\{\omega^1,~\cdots~,\omega^r\},$ then the above condition can be restated by the equality $\omega^1\wedge~\cdots~\wedge\omega^r\wedge\Omega^{h+1}=0$ for the exterior forms defined on open sets contained in the domain of the generators $\omega^i.$ The \textbf{gender} of a Pfaffian system $\mathcal{P}$ is equal to the gender of the family $d\omega$ where $\omega$ is an arbitrary local section of $\mathcal{P}$ and, of course it is given by the gender of any system of generators We also observe that the system $\mathcal{P}$ is integrable if and only if it is of gender zero (\cite{Cartan1901}).

\section{The Cartan co-variant approach}
We start by considering the linear $1-$dimensional elements belonging to $\mathcal{P}$ namely, all the $1-$dimensional sub-spaces of $\mathcal{P}$ at the point $x_0$ and inquire which should be considered as 1-dimensional \textit{integral elements}. Let $\omega$ be a generator of such a linear element and let us imagine that its annihilator $ker~\omega_{x_0}$ is an integral element for the Pfaffian system generated by $\omega.$ Then the condition reads $(d\omega)_{x_0}\equiv 0\hspace{3 mm}mod~\omega_{x_0}.$ It seems however more adequate, for reasons that shall be clarified hereafter, 
to assume the stricter requirement for the integral elements namely,\footnote{Apparently, Cartan missed something (pg.4, \textit{l.}10-12) but, most probably, the present author did not read the hidden lines.}
\begin{equation}
(d\omega)_{x_0}\equiv 0\hspace{3 mm}mod~\mathcal{P}.
\end{equation}
In terms of local generators, the above condition reads
\begin{equation*}
(d\omega)_{x_0}\equiv 0\hspace{3 mm}mod~\{\omega_{1,x_0},~\cdots~,\omega_{r,x_0}\}.
\end{equation*}

\vspace{3 mm}
\noindent
Next, we inquire which should be considered as being the $2-$dimensional integral elements containing a given (fixed) $1-$dimensional element. Though entering in conflict with the notations adopted previously in the \textit{contra-variant} discussion, we shall keep these notations in order to conform to Cartan's writing. We choose a given (fixed) though arbitrary linear integral element $E_1$ and consider the sub-space of all those linear integral elements that are \textit{in involution} with the given element. Inasmuch as in the contra-variant setting, not all the linear integral elements are suitable and the condition for this to be so is, of course, the integrability condition involving the two elements. More precisely, the co-vector $\omega_{x_0},$ or the sub-space generated by it, is in involution with $E_1$ when, by definition,
\begin{equation}
\omega_{x_0}\equiv 0\hspace{3 mm}mod~\mathcal{P},\hspace{3 mm}(d\omega)_{x_0}\equiv 0\hspace{3 mm}mod~\mathcal{P}\hspace{3 mm}and 
\end{equation}
\begin{equation*}
d\omega_0\equiv 0\hspace{3 mm}mod~\mathcal{P}, 
\end{equation*}
where $\omega_0$ is chosen so as to induce a non vanishing co-vector belonging to $E_1.$ In other words, this simply means that $\omega\in\mathcal{P}$ and that
\begin{equation*}
d\omega\wedge\varphi^1\wedge~\cdots~\wedge\varphi^r=0,
\end{equation*}
the same holding for $\omega_0.$ We can now choose $E_2$ containing $E_1,$ only composed of co-vectors that are in involution with $E_1$ and consider the space of all those co-vectors, belonging to $\mathcal{P},$ that are simultaneously in involution with all the elements of $E_2.$ Continuing this process, we choose an element $E_3$ containing $E_2$ only composed by co-vectors that are in involution with all the elements of $E_2$ and so forth. The $character~s_1$ defined by Cartan as being the integer satisfying $n-r-s_1=dim~\tilde{E}_1,$ where the latter is the space of all the integral elements in involution with $E_1,$ is then equal to $\rho_k$ in the contra-variant setting and both situations, the co-variant and the contra-variant, are mirror images one of the other with respect to $\mathcal{P}.$ Furthermore, the result of the previous section concerning the derived system has, in this context, a very simple proof. In the next sections, we discuss Pfaffian systems with characters larger than one that provide a wide spectrum of different situations. We shall nevertheless keep within the Cartan co-variant setting since, as already mentioned earlier, the discussion as well as the calculations become straightforward and simple though we stop understanding anything since we are rowing against the tide.

\vspace{5 mm}
\noindent
From what was shown in the previous section, we can state the following 

\newtheorem{kiss}[LemmaCounter]{Lemma}
\begin{kiss}
A necessary and sufficient condition that the character of a Pfaffian system $\mathcal{P}$ be equal to one is that the rank of its derived system $\mathcal{P}_1$ be one unit less than its own rank.
\end{kiss}

\section{Integration of Pfaffian systems with character equal to one}
We start by taking a local basis $\{\omega^i\}$ for the Pfaffian system $\mathcal{P}$ such that $\{\omega^2,~\cdots~,\omega^r\}$ generates its associated derived system $\mathcal{P}_1,$ complete this basis to a local co-frame by adding the 1-forms $\{\overline{\omega}^j\}$ and continue to assume that the characteristics of $\mathcal{P}$ are null.

\vspace{5 mm}
\noindent
Let us first prove the result that Cartan considered so much surprising, recalling that \textbf{h} denotes the \textbf{gender} of a Pfaffian system.

\newtheorem{blss}[LemmaCounter]{Lemma}
\begin{blss}
A necessary and sufficient condition that the first derived system $\mathcal{P}_1$ of a Pfaffian system $\mathcal{P},$ of character 1, be integrable is that $\textbf{h}\geq 2.$
\end{blss}

\vspace{5 mm}
\noindent
Adapting the generators of $\mathcal{P}$ in such a way that the forms $\{\omega^i\},~i\geq 2,$ generate $\mathcal{P}_1,$ we have by definition
\begin{equation}
d\omega^i\equiv 0\hspace{3 mm}mod~\{\omega^1,~\cdots~,\omega^r\},\hspace{3 mm}i\geq 2,
\end{equation}
and write, to begin,
\begin{equation}
d\omega^2\equiv\chi\wedge\omega^1\hspace{3 mm}mod~\{\omega^2,~\cdots~,\omega^r\}
\end{equation}
where $\chi$ can be taken to be a linear combination of the forms $\overline{\omega}^j.$ A second differentiation then yields
\begin{equation*}
d\chi\wedge\omega^1+\chi\wedge d\omega^1\equiv 0\hspace{3 mm}mod~\{\omega^2,~\cdots~,\omega^r\}    
\end{equation*}
hence
\begin{equation*}
\chi\wedge d\omega^1\equiv 0\hspace{3 mm}mod~\{\omega^1,~\cdots~,\omega^r\}    
\end{equation*}
on account of (7). We next assume that $\chi$ and $\omega^1$ are independent. Since $\chi$ only contain terms in $\overline{\omega}^j,$ we can write $\chi\wedge d\omega^1=\mu+\eta$ where each term of $\mu$ contains some element $\omega^i$ and $\eta$ only contains terms that are expressed by means of the $\overline{\omega}^j.$ Moreover, the term $\eta$ is equal to the product $\chi\wedge\overline{\eta},$ where $\overline{\eta}$ is the sum of all the terms in $d\omega^1$ that are written only with the help of the forms $\overline{\omega}^j.$ It then follows that the exterior product $d\omega^1\wedge d\omega^1$ presents, in its unique term not containing any form $\omega^i,$ a double product $\chi\wedge\chi$ hence this term vanishes and the double product above becomes equal to:
\begin{equation*}
d\omega^1\wedge d\omega^1\equiv 0\hspace{3 mm}mod~\{\omega^1,~\cdots~,\omega^r\}.   
\end{equation*}
Consequently the gender of $\mathcal{P}$ is at most equal to 1 contradicting the initial hypothesis. We infer that $\chi$ and $\omega^1$ are dependent and that the relation (8) reduces consequently to
\begin{equation}
d\omega^2\equiv 0\hspace{3 mm}mod~\{\omega^2,~\cdots~,\omega^r\}.
\end{equation}
The same argument being valid for all the other forms generating $\mathcal{P}_1,$ the integrability of this derived system then follows. 

\vspace{5 mm}
\noindent
We can now proceed with the integration of the system $\mathcal{P}$ namely, construct integral manifolds of dimension $n-r-1$ that are tangent to given maximal integral contact elements, also of the samee dimension, since the character of $\mathcal{P}$ is assumed to be equal to one. Fixing the point $x_0,$ we denote by $\mathcal{L}$ the integral leaf of $\mathcal{P}_1$ that contains the given point. Then $dim~\mathcal{L}=n-r+1$ and, if we assume as previously that the forms  $\omega^i,~i\geq 2,$ generate $\mathcal{P}_1,$ the restricted system $\mathcal{P}|\mathcal{L}$ is generated by the single 1-form $\omega^1|\mathcal{L}.$ Obviously, any maximal integral element issued at the point $x_0$ is contained in $T_{x_0}\mathcal{L}$ and any maximal integral manifold of $\mathcal{P}$ containing $x_0$ is contained in $\mathcal{L}$ being as well an integral manifold of the restricted system. Moreover, the maximality properties are preserved under restriction since the restricted system is not integrable, otherwise there would be an $(n-r)-$dimensional integral manifold of $\mathcal{P}$ passing through $x_0$. Since the \textit{gender} of a Pfaffian system is a point-wise notion involving only the exterior algebra of the co-tangent spaces at given points, we infer that the gender of the restricted system is the same as that of the initial system. On the other hand, the restricted system $\mathcal{P}|\mathcal{L}$ is non-integrable and has rank equal to one hence is a Darboux system (\cite{Kumpera2014}). The definition \textit{per se} of the gender of this restricted system will then assert that $\mathcal{P}|\mathcal{L}$ is a Darboux system of Cartan or, inasmuch, Darboux class equal to  $2\textbf{h}+1$ hence generated locally by a form
\begin{equation*}
dz^{h+1}+p^1dz^1+p^2dz^2+~\cdots~+p^hdz^h,    
\end{equation*}
where $\it{h}=\textbf{h}$ is the gender at a point\footnote{Cartan preferred minus signs.}. Putting all the stacks together, we infer that $\mathcal{P}$ admits local models of the form:
\begin{equation*}
dy^1\hspace{64 mm}
\end{equation*}
\begin{equation*}
dy^2\hspace{64 mm}
\end{equation*}
\begin{equation}
..........\hspace{63 mm}
\end{equation}
\begin{equation*}
dy^{r-1}\hspace{59 mm}
\end{equation*}
\begin{equation*}
dz^{h+1}+p^1dz^1+p^2dz^2+~\cdots~+p^hdz^h
\end{equation*}
with respect to suitable local coordinates and where the functions $y^1,y^2,~\cdots~,y^{r-1}$ are chosen to be independent first integrals of the (integrable) derived system (\textit{c.f.}, \cite{Cartan1901}, pg.28). The solutions of the above system will provide, locally, the maximal integral manifolds of $\mathcal{P}.$ If by any chance the given Pfaffian system admits non-vanishing characteristics\footnote{Cartan never imposed any restriction on the characteristics since, in his \textit{modus operandi}, nothing much would change.}, then it suffices to chose, for $y, z$ and $p,$ first integrals of the characteristic system.

\vspace{5 mm}
\noindent
Let us finally examine what happens when $\textbf{h}=1.$ We can in this case integrate, starting from a given point $x_0$ on $\mathcal{L}$, all the vector fields belonging to the annihilator of $\mathcal{P}|\mathcal{L}.$ Denoting by 
\begin{equation*}
t\longmapsto(z^1(t),~\cdots~,z^{\mu}(t)),\hspace{5 mm}\mu=n-r+1,
\end{equation*}
a generic integral curve commencing at $x_0,$ we plot all the points obtained by considering the images of the above integral curves conditioned to a given specific, though arbitrary, non-trivial relation among the targets $z^i.$ A detailed proof of this statement is rather long and shall be omitted. Fortunately, the author is backed up by Cartan's fantastic intuition (\cite{Cartan1901}, pg.28-29, where he refers to von Weber \cite{Weber1898}).

\section{Pfaffian systems whose characters are larger than one}
We begin by examining systems with characters equal to 2 and, to simplify, we assume right away that their characteristics are null. Much will then depend upon the \textit{gender} of the system and offers several options. We exhibit firstly conditions under which the character being equal to 2 implies that the rank of the first derived system is equal to $r-2,$ the general rule not being as strict in the present case as it was previously. It seems furthermore worthwhile, at this point, to make a few comments since our statement is, apparently, in disagreement with Cartan's claims. Assume, for a moment, that the Pfaffian system $\mathcal{P}$ has non-trivial characteristics and denote by $\mu$ the difference in the dimensions of the initially given manifold \textit{M} and the (local) quotient manifold $\overline{M},$ modulo the characteristics of the system $\mathcal{P}$. Then $\mathcal{P}$ factors to an equivalent system $\overline{\mathcal{P}},$ both systems have the same rank and, consequently, the annihilator distributions $\Sigma$ and $\overline{\Sigma}$ differ, in their point-wise dimensions, by $\mu,$ the dimension of $\overline{\Sigma}$ being lesser since the characteristic distribution is contained in $\Sigma$. Similar statements hold for the derived system and its quotient in $\overline{M}$ though the latter can still have non-trivial characteristics thus explaining the discrepancy concerning the rank of $\overline{\mathcal{P}_1}$ and the dimensions of the maximal integral elements of $\mathcal{P}$ and $\overline{\mathcal{P}}.$ We also observe that the passage to the derived systems is functorial and therefore compatible with quotients \textit{i.e.}, $\overline{\mathcal{P}_1}=(\overline{\mathcal{P}})_1.$

\vspace{5 mm}
\noindent
In the sequel we shall only give a brief account on Cartan's results so as to promptly continue with our main discussion. In what concerns the Pfaffian systems of character 2, the main result is the following:

\newtheorem{buss}[TheoremCounter]{Theorem (Cartan)}
\begin{buss}
The rank of the derived system of a Pfaffian system $\mathcal{P}$ of character two is always two units less than the rank of the system except when $\mathcal{P}$ admits characteristic elements with dimensions not less than  $dim~\Sigma_{x_0}-3$ in which case the rank lowers by at least three units.
\end{buss}

\noindent
The proof is similar to that of the previous lemma as soon as we choose two independent local sections of the system whose differentials do not vanish on the kernel and an entirely analogous result also holds in general. We shall therefore only give some attention to the second part of the statement that is far more delicate. Let us then assume, the Pfaffian system being of character 2, that the rank of the derived system $\mathcal{P}_1$ is equal to $r-\mu$ where $\mu\geq 3.$ This being so, we can choose a local basis $\{\omega^1,~\cdots~,\omega^{\mu},~\cdots~,\omega^r\}$ for $\mathcal{P}$ and in a neighborhood of the point $x_0$ in such a way that the forms $\{\omega^1,~\cdots~,\omega^{\mu}\}$ do not belong to $\mathcal{P}_1,$ hence inasmuch any non-trivial linear combination of these forms cannot belong to the derived system, the remaining forms generating $\mathcal{P}_1.$ This being so and since $d\omega^j\not\equiv 0\hspace{3 mm}mod~\mathcal{P},~j\leq\mu,$ we can select $\mu$ independent vectors $v_j$ in $\Sigma_{x_0}$ not belonging, of course, to the characteristic sub-space of $\mathcal{P}$ at the point $x_0$ in such a way that the family
\begin{equation*}
\{\omega^1,~\cdots~,\omega^r,i(v_1)d\omega^1,~\cdots~,i(v_{\mu})d\omega^{\mu}\}
\end{equation*}
is a local basis for the characteristic system of $\mathcal{P}.$ The rank $r+\mu$ of this characteristic system must be equal to \textit{n}, hence $\mu =n-r,$ otherwise the system $\mathcal{P}$ would have non-trivial  characteristics contradicting our assumption. But then, the system $\mathcal{P}$ will be of rank equal to the dimension of \textit{M}, its annihilator is the null distribution and the integral leaves are the points of \textit{M}, a rather uninteresting foliation. Within the scope of our assumptions, we can therefore restate the above theorem as follows:

\newtheorem{boss}[TheoremCounter]{Theorem}
\begin{boss}
Let $\mathcal{P}$ be a Pfaffian system with null characteristics and character equal to 2. Then the rank of its  derived system is two units less than the rank of the system.
\end{boss}

\noindent
Let us now say a word on the number 3. Apart from being an extraordinary prime number that brings much luck, let us argue as follows:

\noindent
We select among all the possible non-trivial linear combinations of the forms $\{\omega^1,~\cdots~,\omega^{\mu}\}$ one such that has a minimum gender at the point $x_0.$ Re-arranging once more the local basis, we can assume for convenience that the above form is simply $\omega^1$ and can therefore write, at the point $x_0,$ $(d\omega^1)^{\textbf{h}+1}\equiv 0\hspace{3 mm}mod~\{\omega^1,~\cdots~,\omega^r\},$ $\textbf{h}$ being the minimum value for which the above congruence holds. A simple calculation will then show that the value $\textbf{h}=1$ is the only possible value for which the rank of the characteristic system is  larger than $r-2$ and consequently the rank of the derived system is lesser than $r-2,$ precisely equal to $r-3$. 

\vspace{5 mm}
\noindent
However, we owe presently our apologies to Élie Cartan for not considering systems with non vanishing  characteristics. It so happens that in many applications Pfaffian systems with non-trivial characteristics do appear most naturally and, furthermore, the factoring of a Pfaffian system, modulo its characteristics, is never an easy task since it requires an integration process that only reduced locally to a system of ordinary differential equations when the characteristic leaves are $1-$dimensional.

\vspace{5 mm}
\noindent
In the sections 29-34, Cartan discusses those Pfaffian systems for which half of the point-wise dimension of their annihilators $\Sigma$ is not less than the value of their characters. Since, once more, the discussion of this setup involves non-vanishing characteristics, we shall leave it aside. There is however a special case when the equality of the above mentioned values holds and that merits to be commented. Let us then assume that the previously mentioned system $\mathcal{P}$ also verifies the following properties:

\vspace{5 mm}
(a) The character of $\mathcal{P}$ is equal to 2,
 
\vspace{3 mm}
(b) $\mathcal{P}$ has null characteristics and 
 
\vspace{3 mm}
(c) All the maximal integral contact elements that contain a specific $1-$dimensional integral element $E_1$ also contain an integral element $E_2\supset E_1$ whose dimension is not less than 2.

\vspace{5 mm}
\noindent
Under these conditions, the equality $\frac{n-r}{2}=\rho_k$ holds. Translated in terms of maximal integral manifolds, the condition (c) means that if two such manifolds have in common a line they also have in common at least a surface. Cartan calls such systems \textit{systatical (systatiques)}.

\section{Singular Pfaffian Systems}
For generic Pfaffian systems, the successive characters assume their maximum values hence, in particular, for systems with character equal to 2 the values are:
\begin{equation*}
s_1~=~s_2~=~\cdots~=~s_{\rho_k-1}~=~2.    
\end{equation*}
\textit{Quant à} $s_{\rho_k}$ \textit{il est égal a zéro, un ou deux, suivant le reste de la division de $n-r$ par 3 (Cartan dixit)}.

\vspace{5 mm}
\noindent
This being so, $s_n=0$ when $n-r=3\rho_k-2,$ $s_n=1$ when $n-r=3\rho_k-1,$ and $s_n=2$ when $n-r=3\rho_k.$ 

\vspace{5 mm}
\noindent
We shall say, together with Cartan, that the Pfaffian system $\mathcal{P},$ with character equal to 2, is \textit{singular} when $s_{\rho_k-1}\leq 1.$ Still quoting Cartan, for a singular system having character two and null characteristics the maximum dimension of the linear integral contact elements \textit{i.e.}, the value of $\rho_k$ is at least equal to 3 and, consequently, the point-wise dimension of the distribution $\Sigma$ is at least equal to 6. Another rather surprising property, consequence of the above assertions, is the following:

\newtheorem{bosta}[PropositionCounter]{Proposition}
\begin{bosta}
All the systatical Pfaffian systems are singular as soon as $\rho_k\geq 3.$
\end{bosta}

Cartan's \textit{Mémoire} then terminates with an attempt to extend some of the above properties to non-systatical systems. A very long and typically end of the $19^{th}$ beginning of the $20^{th}$ century discussion culminates in the following statement where "\textit{can be integrated}" means the obtainment of the maximal integral manifolds (pg.70):

\vspace{5 mm}
\noindent
\textit{A singular Pfaffian system with null Cauchy characteristics and character equal to 2 can be integrated by means of a system of ordinary differential equations. In case the system admits non-trivial characteristics, its characteristic system can also be integrated by means of ordinary differential equations whenever the derived system of the given system is non-integrable.}

\vspace{5 mm}
\noindent
Curiously enough, the second order \textit{Monge characteristics} play an important role in the proofs.

\vspace{5 mm}
\noindent
Finally we arrive at what really matters namely, the integration of the above considered Pfaffian systems via the reductions provided by Jordan-Hölder sequences. In the previous discussion, we have seen to what extent the derived systems are relevant. Consequently, we begin by integrating the derived system and afterwards locate, among the manifolds thus obtained, those that are maximal for the given system. We are thus faced in constructing Jordan-Hölder resolutions for the derived system and each specific situation will exhibit its particular techniques. It might however be of some advantage to try initially integrating one of the systems that make part of the sequence of iterated derived systems and, in this case, the chain (3) will already provide some useful terms. In particular, it is possible to detect the appearance of ordinary differential equations, 
these showing up when the quotients are $1-$dimensional. \nocite{Kumpera1999}

\section{Examples}
Allowing initially $dim~M=5$ and $rank~\mathcal{P}=3,$ the system generated by
\begin{equation*}
\{\omega^1=dx^1+x^4dx^5,~\omega^2=dx^2,~\omega^3=dx^3\}
\end{equation*}
is \textit{not} integrable since
\begin{equation*}
d\omega^1=dx^4\wedge dx^5\not\equiv 0\hspace{3 mm}mod~\{\omega^1,\omega^2,\omega^3\}
\end{equation*}
and $\mathcal{P}_1$ is generated by $\{\omega^2,\omega^3\}$ hence is, of course, integrable. The character of $\mathcal{P}$ is equal to 1 since, at any point $x_0,$ the maximal integral manifold is the line defined by 
$x^1=c^1,x^2=c^2,x^3=c^3,x^5=c^5$ for some constants $c^i$ and is a line parallel to the $x^4-$axis. It should be observed that the annihilator $\Sigma_{x_0}$ of $\mathcal{P}_{x_0}$ is $2-$dimensional.

\vspace{5 mm}
\noindent
Let us now replace $\omega^3$ by $\tilde{\omega}^3=dx^3+x^5dx^1.$ Then
\begin{equation*}
d\tilde{\omega}^3=dx^5\wedge dx^1=dx^5\wedge\omega^1,
\end{equation*}
hence $\mathcal{P}_1$ is again generated by $\{\omega^2,\omega^3\}$ but is \textit{not} integrable since $d\tilde{\omega}^3\not\equiv 0\hspace{3 mm}mod~\{\omega^2,\omega^3\}.$ This derived system, whose annihilator is point-wise $3-$dimensional, just admits the $2-$dimensional integral manifolds defined by the equations $x^1=c^1,x^2=c^2,x^3=c^3.$

\vspace{5 mm}
\noindent
We next give a glance at those Pfaffian systems with character equal to 2 and it will just suffice to add one more dimension to our space. Let us therefore consider the system $\mathcal{P}$ generated, on a $6-$space, by
\begin{equation*}
\{\omega^1=dx^1+x^4dx^5,~\omega^2=dx^2+x^5dx^6,~\omega^3=dx^3\}    
\end{equation*}
Then the annihilator $\Sigma_{x_0},$ of $\mathcal{P},$ is $3-$dimensional, the derived system $\mathcal{P}_1$ is generated by $\{dx^3\}$ hence is integrable and its rank is two units less than the rank of $\mathcal{P}.$ Consequently, $\mathcal{P}$ is not integrable. On the other hand, the generators vanish on the affine line defined by the equations
\begin{equation*}
x^1=c^1,~x^2=c^2,~x^3=c^3,~x^5=c^5,~~x^6=c^6
\end{equation*}
and, moreover, this line is a maximal integral manifold. The character of $\mathcal{P}$ is therefore equal to 2.

\vspace{5 mm}
\noindent
Let us now replace $\omega^3$ by $\tilde{\omega}^3=dx^3+x^5dx^1.$ Then $\mathcal{P}$ continues to be of character 2 but the derived system $\mathcal{P}_1,$ generated by $\tilde{\omega}^3,$ is no longer integrable since $d\tilde{\omega}^3=dx^4\wedge dx^1\not\equiv 0\hspace{3 mm}mod~\{\tilde{\omega}^3\}$ and only admits the $4-$dimensional maximal integral manifolds ($4-$dimensional affine spaces) defined by the equations
\begin{equation*}
x^1=c^1,~x^3=c^3.
\end{equation*}
The maximal integral element of $\mathcal{P}$ are only $1-$dimensional and its maximal integral manifolds are affine lines.

\vspace{5 mm}
\noindent
We continue our little game by constructing now a Pfaffian system of character 2 whose derived system has rank three units less than the rank of the given system \textit{i.e.}, we shall construct a singular Pfaffian system with character equal to 2. For this, we simply consider the system $\mathcal{P}$ generated, on a $6-$dimensional space, by
\begin{equation*}
\{\omega^1=dx^1+x^4dx^5,~\omega^2=dx^2+x^5dx^6,~\omega^3=dx^3+x^6dx^4\}.
\end{equation*}
The maximal integral manifolds of $\mathcal{P}$ are affine lines and a direct calculation shows that $\mathcal{P}_1=0.$ The rank of the derived system is three units less than its own rank.

\vspace{5 mm}
\noindent
As for the \textit{gender} of a Pfaffian system, the first example considered above illustrates the content of the Lemma 4 since the local section $\omega^1+x^2\omega^3$ has its gender equal to 2.

\vspace{5 mm}
\noindent
Let us finally terminate our discussion by exhibiting three \textit{Jordan-Hölder resolutions}. To begin with, we continue to consider the first example in $5-$space and observe that the pseudo-group $\Gamma$ of all the local automorphisms of $\mathcal{P}$ leaves also invariant the derived system $\mathcal{P}_1$ hence, inasmuch, its annihilator $\Sigma_1$ whose point-wise dimension is equal to 3. Needless to say that the integral manifold of $\Sigma_1$ are all parallel affine three dimensional sub-spaces in $5-$space. Moreover, $\Gamma$ also leaves invariant the annihilator $\Sigma\subset\Sigma_1.$ We denote by $\Gamma_1$ the sub-pseudogroup of all the elements $\varphi,$ belonging to $\Gamma,$ that \textit{maintain invariant} each connected integral manifold of $\Sigma_1$ \textit{i.e.}, $\varphi$ transforms the points of such an integral manifold into points of the same manifold. Equivalently, this simply means that the above mentioned integral manifolds are the intransitivity classes of $\Gamma_1.$ We can now factor, locally, the given $5-$ modulo these intransitivity classes and obtain, as quotient, an open set \textit{U} in numerical $2-$space. Furthermore, the Pfaffian system $\mathcal{P}$ being invariant under $\Gamma_1,$ factors to the quotient and yields a quotient Pfaffian system $\overline{P}$ that is generated by $\{dx^1,~dx^2\}$ hence its maximal integral leaves are the points of \textit{U}. Fixing a point $y_0\in U$ and denoting by $Y_0$ the inverse image $q^{-1}(y_0),$ where \textit{q} is the quotient projection modulo the above intransitivity classes, we can now restrict the system $\mathcal{P}$ to the sub-manifold $Y_0,$ this restricted system $\mathcal{Q}$ being generated by the $1-$form $dx^1+x^4dx^5.$ Since $dim~Y_0=3$ and since the generating form of the system is of maximum Darboux class equal to 3, the restricted system is a Darboux system of rank 1 and only admits $1-$dimensional integral manifolds that, in fact, are the affine lines in $5-$space that constitute the maximal integral manifolds of $\mathcal{P}.$ The integration of this system is consequently achieved via the Jordan-Hölder method since the quotient pseudo-group $\Gamma/\Gamma_1$ is equal to the set of all the local transformations in \textit{U}, hence is simple.

\vspace{5 mm}
\noindent
We shall continue with the same Pfaffian system but apply the Jordan-Hölder methodology in a different manner. The character of the system being equal to 1 and, consequently, the maximal integral manifold being $1-$dimensional curves that actually constitute a foliation, we consider the sub-pseudogroup $\Gamma_1$ as being composed by all those local transformations that keep invariant the connected maximal integral curves, in much the same way as considered previously for the integrals of the first derived system. As before, let \textit{U} be a local quotient of the $5-$space modulo the above integral curves. Then the dimension of \textit{U} is equal to 4 and we can, inasmuch, factor locally the system $\mathcal{P}$ to the space \textit{U}. Since the annihilator, at each point, of the system $\mathcal{P}$ contains the tangent space to the integral curve at that point, we infer that the quotient system $\mathcal{Q}$ is still of rank 3 and is generated by $\{dx^2,~dx^3\}$ together with the quotient of $\omega^1.$ However, since the dimension of \textit{U} is equal to 4 and the quotient system $\mathcal{Q}$ has rank equal to 3, we infer that it is necessarily integrable, its maximal integrals being again curves. Proceeding as prescribed in the Jordan-Hölder scheme, we fix a certain, though arbitrary, integral curve of $\mathcal{Q}$ and consider the inverse image of this curve (more precisely, of the image of this curve) modulo the established quotient projection described above. In other terms, we consider the surface \textit{Y} obtained as the union of all the integral lines of $\mathcal{P}$ that project onto points of the selected (though fixed) integral curve of $\mathcal{Q}.$ Knowing \textit{a priori} that the integral curves are parallel affine lines, we shall thus obtain a \textit{ruled surface}. Let us now restrict $\mathcal{P}$ to the above surface and denote this restricted Pfaffian system by $\overline{\mathcal{P}}.$ Since it must vanish on the tangent spaces to the integral curves, its rank must be equal to 1 and its integral curves are precisely the integral curves of $\mathcal{P}$ contained in the above surface. The Jordan-Hölder method is therefore a most outstanding integration procedure since it enables us to obtain the maximal integral manifolds of the initially given system by a very precise and constructive technique. It should be noted, inasmuch, that the whole process can of course be performed without knowing \textit{a priori} the nature of the maximal integral manifolds that we are looking for (curves in our present context), this being most often the case and  in fact the sole purpose of an integration method. Moreover, for the system discussed here, the usual  \textit{quadrature} technique (Calculus I), will enable us to find a first integral for the integral curves.  

\vspace{5 mm}
\noindent
Let us now consider the last example given above, where the three forms $\omega^i,$ generating the system, are of Darboux class equal to 3. We proceed as above and consider the sub-pseudogroup $\Gamma_1$ composed by all  those elements of $\Gamma$ that preserve the maximal integral manifolds of $\mathcal{P},$ here again $1-$dimensional curves since the character of $\mathcal{P}$ is equal to 2. The local quotient space \textit{U} is, in this case, $5-$dimensional and the resulting quotient system $\mathcal{Q}$ continues to be of rank 3, the point-wise dimension of its annihilator being thereafter equal to 2. However, this quotient system acquires  character equal to 1.\footnote{somehow, the reduction in the dimension squeezes the system into a smaller ambient space.} Taking a $1-$dimensional integral curve $\gamma$ of $\mathcal{Q},$ considering the inverse image $Y=q^{-1}(im~\gamma)$ defined in the initially given $6-$space and restricting $\mathcal{P}$ to the surface \textit{Y}, we finally obtain a $rank~1$ Pfaffian system whose integral curves are precisely the integral curves of $\mathcal{P}$ contained in \textit{Y} and for which a first integral can be obtained by a quadrature.

\vspace{5 mm}
\noindent
We should at present confess to the reader that the above discussion is actually \textit{a big fake}. In fact, the aim of the Jordan-Hölder procedure is, of course, the obtainment of maximal integral manifolds of the given system whereas our procedure was to exhibit this technique, on examples, with the help of the \textit{a priori} knowledge of the solutions. Nevertheless, we believe that the reader is able to devise much more sophisticated examples where the true aim will be, of course, the determination of maximal integral manifolds not known in advance. One last word is due. Élie Cartan never claimed that the Jordan-Hölder procedure was an easy matter. It can in fact become extremely involved, reason for which Cartan \textit{never} mentioned neither Jordan nor Hölder.

\vspace{5 mm}
\noindent
A very interesting type of non-integrable Pfaffian systems is given by the \textit{Flag Systems} considered in various occasions by Élie Cartan (\textit{e.g.}, \cite{Cartan1901}, \cite{Cartan1910}, \cite{Cartan1911}). We claim that a flag system $\mathcal{F}$ with null characteristics has always its character equal to $n-r-1,$ where $n=dim~M$ and $r=rank~\mathcal{F}.$ In fact, proceeding as indicated in \cite{Kumpera2014} and replacing the descending chain of \textit{derived} systems by the more appropriate \textit{structured multi-fibration} (6) where each successive derived system is replaced by an equivalent system that now has also null characteristics, we see that the maximal dimensional integral manifolds, at each stage, are the integral curves contained in the inverse image of the integral curves of the Darboux system, these last curves being also the maximal dimensional integral curves. Our claim then follows since the point-wise dimension of the annihilator distribution $\mathcal{F}^{\perp}$ is precisely $n-r.$ More generally, we can consider the \textit{Multi-Flag Systems} or, still better, the \textit{Truncated Multi-Flag Systems} that have much relevance in the study of under-determined ordinary differential equations in what concerns the \textit{Monge property} (\cite{Kumpera2002}). We leave the details to the reader and, in particular, the calculation of the character for such systems. It should be noted that the difference between flag and multi-flag systems can be recognized by looking at their local equivalence pseudo-groups that are non-isomorphic. Whereas for flag systems the Monge condition is a necessary and sufficient condition for the validity of the Monge property (\cite{Cartan1911}), for multi-flag systems the condition is only sufficient (\cite{Kumpera2002}). 

\vspace{5 mm}
\noindent 
Last but not least, let us give a glance at a rare jewel left to us by Élie Cartan namely, his article on Galois theory (\cite{Cartan1938}), where he bases his argumentation on the classical Picard-Vessiot theory for linear differential systems and shows how to integrate such systems via the Jordan-Hölder method.

\vspace{5 mm}
\noindent
On th base manifold \textit{M}, we consider the $k-$th order jet space $J_kTM$ of all the $k-$jets of local sections of the tangent bundle $\pi:TM~\longrightarrow~M,$ \textit{i.e.}, $k-$jets of local vector fields on \textit{M}. Then $\alpha_k:J_kTM~\longrightarrow~M,$ $\alpha_k$ denoting the \textit{source} map, becomes again a vector bundle, its linear operations being simply the extensions, to $k-$jets, of the linear operations in \textit{TM}. A linear partial differential equation of order \textit{k} is, by definition, a vector sub-bundle $\mathcal{R}$ of $J_kTM.$ The dual (covariant) version reads as follows. We take the dual bundle $T^*M,$ consider the jet space $J_kT^*M$ and then put in evidence a vector sub- bundle $\mathcal{S}.$ It should be remarked that $J_kT^*M\equiv(J_kTM)^*.$ Given a linear differential equation $\mathcal{R}\subset J_kTM,$ its annihilator becomes a vector sub-bundle of $J_kT^*M$ \textit{i.e.}, $\mathcal{R}^{\perp}=\mathcal{S}\subset J_kT^*M$ and, inasmuch, we can consider the annihilator, in $J_kTM,$ of any sub-bundle $\mathcal{S}\subset J_kT^*M.$ We denote $n=dim~M$ and $n_k=dim~J_kTM.$ On account of the previous remarks, it also follows that $n_k=dim~J_kT^*M.$ Taking a local coordinate system $(U,(x^i))$ in the base space \textit{M}, we derive corresponding coordinate systems for $J_kTM$ and $J_kT^*M,$ defined on the respective inverse images of the open set \textit{U}, and any linear partial differential equation will be determined, locally by a system of linear equations with variable coefficients. If $rank~\mathcal{R}=\rho$ denotes the point-wise dimension of the fibres, then $rank~\mathcal{S}=n_k-\rho$ and the number of independent linear equations defining locally the system is also equal to $n_k-\rho.$

\vspace{5 mm}
\noindent
As is usual, we indicate by $\mathcal{C}_k$ the canonical contact system defined on $J_kTM$ and, just to mark the distinction, by $\mathcal{C}_k^*$ the canonical contact system defined on $J_kT^*M.$ We recall that a local section in the $k-$jets is holonomic \textit{i.e.}, results from a section in the base tangent or co-tangent bundles if and only if it is annihilated by the contact Pfaffian syatem and, furthermore, a local transformation in the jet space is the prolongation of a base space transformation (in $TM$ or $T^*M$) if and only if it preserves the corresponding contact structure. In the present case, we shall restrict our attention to the pseudo-group of all the local transformations that are vector bundle morphisms. This is the starting point of the Picard-Vessiot theory and the set of all the vector bundle automorphisms of the equation is called its \textit{Galois group}. Most often, this pseudo-group operates non-transitively on the equation though the $(k+1)-st$ order  groupoid of all the $(k+1)-$jets that preserve the equation, via the \textit{semi-holonomic} action, is always transitive on account of linearity. This action is defined as follows. Given a jet $j_{k+1}\varphi(v)$ on $TM$ or $T^*M,$ we consider the flow $j_k\varphi:w~\longmapsto~j_k\varphi(w)$ and subsequently take the jet $j_1(j_k\varphi)(X),~\alpha_k(X)=v.$ Restricting the contact structures to the corresponding equations, we obtain the Pfaffian systems $\overline{\mathcal{C}}_k$ and $\overline{\mathcal{C}}_k^*$ associated to the linear equations. These systems are not integrable but their $n-$th dimensional integral manifolds transversal to the $\alpha_k-$ fibres are precidely the images of the $k-$jets of solutions. On account of the transitivity at the groupoid level, we can apply the Jordan-Hölder integration method as outlined in \cite{Kumpera2016}. When the linear pseudo-group operates transitively on the equations, the integration process can be achieved by just calling out the simple linear Lie groups. We finally observe that Picard (\cite{Picard1896}), Vessiot (\cite{Vessiot1904},\cite{Vessiot1912}) and also Drach (\cite{Drach1898}) only considered first order linear systems. Many years earlier, Sophus Lie studied linear differential equations invariant under linear groups with constant coefficients (\cite{Lie1876})  As for Cartan's \textit{mémoire}, he shows us many things, in \cite{Cartan1938}, via examples.

\vspace{5 mm}
\begin{figure}[h!]
\centering
\includegraphics[scale=3.4]{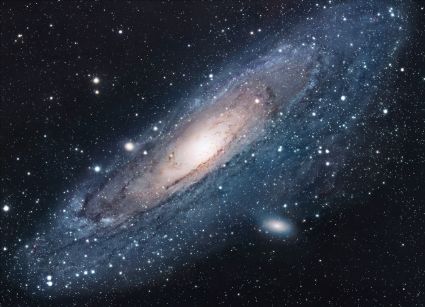}
\caption{The Universe}
\label{fig:univerise}
\end{figure}

\bibliographystyle{plain}
\bibliography{references}

\end{document}